\documentclass[a4paper,11pt]{article}
\usepackage[a4paper,left=2.75truecm,top=3.0truecm,height=22.7truecm,%
width=15.5truecm]{geometry}
\usepackage[T1]{fontenc}
\usepackage{times}
\usepackage{amsmath}
\usepackage{amsfonts}
\usepackage{graphicx}
\usepackage{color}

\newcommand{\0}{{\boldsymbol 0}}
\newcommand{\Alpha}{{\boldsymbol\alpha}}
\newcommand{\eps}{{\varepsilon}}

\newcommand{\PD}[2]{\frac{\partial{#1}}{\partial{#2}}}
\newcommand{\R}{{\Bbb R}}
\newcommand{\PP}{{\mathbb P}}
\newcommand{\Wad}{{W^{ad}}}
\newcommand{\Uad}{{U^{ad}}}
\newcommand{\XX}{{\boldsymbol X}}
\newcommand{\WW}{{\boldsymbol W}}
\newcommand{\calJ}{\mathcal{J}}

\newcommand{\calO}{{\mathcal O}}
\newcommand{\calP}{{\mathcal P}}
\newcommand{\calT}{{\mathcal T}}
\newcommand{\calU}{{\mathcal U}}

\newcommand{\hatOmega}{{\widehat\Omega}}

\newcommand{\nvec}{\boldsymbol n}
\newcommand{\qq}{{\boldsymbol q}}
\newcommand{\rr}{{\boldsymbol r}}
\newcommand{\vv}{{\boldsymbol v}}
\newcommand{\ww}{{\boldsymbol w}}
\newcommand{\xx}{{\boldsymbol x}}
\newcommand{\bomega}{{\overline{\omega}}}
\newcommand{\calOt}{{\widetilde{\cal O}}}

\newcommand{\bld}[1]{\boldsymbol{#1}}

\newtheorem{remark}{Remark}
\newtheorem{example}{Example}

\begin{document}

\title{Topology optimization in Bernoulli free boundary
problems\footnote{The research of the second author was supported by
  the grant \#\,122932 of the Academy of Finland. The third author
  acknowledges the support of the grant IAA100750802 
    of the Grant Agency of the Czech Academy of Science and 
MSM0021620839.
}}

\author{J. I. Toivanen\footnote{Department of Mathematical Information
    Technology, 
  University of Jyv\"askyl\"a, PO Box 35 (Agora), 40014 Jyv\"askyl\"a,
  Finland} \quad
R.A.E. M{\"a}kinen\footnote{Department of Mathematical Information Technology,
  University of Jyv\"askyl\"a, PO Box 35 (Agora), 40014 Jyv\"askyl\"a,
  Finland} 
\quad 
J. Haslinger\footnote{Faculty of Mathematics and Physics,
  Department of  Numerical  Mathematics, Charles 
  University Prague, Sokolovska 83, 186 75, Prague  8, Czech Republic }
}

\maketitle

\begin{abstract}
In this work we consider topology optimization of systems, which
are governed by the external Bernoulli free boundary problem.
We utilize the so-called pseudo-solid approach to solve the governing
free boundary problems during the optimization. To define  design domains
we utilize a level set representation parameterized by radial basis functions.
This design parametrization allows topological changes in the design domain.
\end{abstract}

\section*{Introduction}

In shape optimization, state problems are usually given by partial
differential equations which are formulated and solved in known
domains. Thus the unknown shape we are looking for is   represented by
a solution to our optimization problem. If the state is defined by a
free boundary problem, the situation becomes more involved. The main
feature of such  problems is the fact that not only a function solving
the respective PDE but also the domain itself where this PDE is
considered   are unknown and they have to be found
simultaneously. Unlike to the standard case, the unknown geometry now
appears   in both, the upper as well as the inner level of the
optimization problem. Classical shape optimization with external
Bernoulli free boundary problems (BFBP) as the state relation has been
presented in \cite{toi08}.  A couple $(u,\Omega)$ is said to be a
solution to BFBP if 
$u$ is  harmonic  in a doubly connected (unknown) domain $\Omega$, it
satisfies 
the Dirichlet condition on the  boundary $\partial\omega$ of the inner
component $\omega$ 
given a-priori and, in addition the over-determined system consisting
of the Neumann and Dirichlet condition  on the (unknown) free
boundary. The shape of the free boundary was controlled by $\omega$ but
still keeping the same topology of $\omega$. To make numerical realization
simpler we supposed that the system of admissible $\omega$:s consists of
star-like domains enabling us to express $\partial\omega$ and $\omega$
in terms of the 
polar coordinates.  However it was observed that this system is very
narrow and many reasonable target free boundaries can not be
matched. Indeed, if the boundary $\partial\omega$ is twice
differentiable and $\omega$ is
star-like with respect to a neighborhood of some point in $\omega$ then the
respective free boundary is of the class $C^\infty$
\cite{flu97}. Hence if a given target 
free boundary has the discontinuous curvature (a square with rounded
corners, e.g) then it can never be realized for any such $\omega$. In
computations this fact manifests itself by oscillations of
$\partial\omega$: using  
finer discretizations of $\partial\omega$, the oscillations become faster and faster
(see \cite{toi08}). This behavior indicates the tendency to change the
topology of $\omega$. There are two ways how to overcome such oscillations:
i) to restrict the design space ii) to extend it in such a way that
changes of topology are possible. If i) is used then oscillations are
suppressed but usually there is a big gap between the target and the
found free boundary. To get better results one has to change the
topology of $\omega$.  One of possible ways how to do that will be described
in this paper. We shall use a level set approach. The topology of $\omega$
will determined by the zero level set of a function whose argument is given by
the linear combination of a system of radial basis functions with
overlapping supports. Instead of solving Hamilton-Jacobi equation
describing the evolution of the level set function, the problem can be
treated as a parametric minimization problem with parameters
represented by the coefficients of the above mentioned linear
combination.  

The paper is organized as follows: Section 1 presents a general
setting of a class of topology optimization problems governed by
BFBP. Section 2 is devoted to the presentation of the state solver.
Since the state problem will be solved several times, one has to have
at his disposal an efficient and reliable method. It turns out that
the so-called pseudo-solid approach enjoys both these requirements. In
Section 3 we shortly recall a level set approach widely used in
topology optimization. Section 4 deals with a discretization of the
whole optimization problem and its numerical realization when
$C^2$ radial basis functions are used to
parameterize the level set function. Section~\ref{sec:examples}
presents results of 
several model examples. Finally, the paper is completed with two
appendices on smoothing the Heaviside function and on the analytical 
solution to a specific Bernoulli problem.

Throughout the paper we use the following notation: the symbol
$H^k(\Omega)$ ($k\ge0$ integer) stands for the Sobolev space of
functions which are 
together with their derivatives up to order $k$ square integrable in $\Omega$,
i.e.\ elements of $L^2(\Omega)$ (we set  $H^0(\Omega)\equiv L^2(\Omega)$).

\section{Setting of the problem}

We start with the definition of the state problem represented by an
exterior Bernoulli free boundary problem.
Let $\gamma<0$ and an open set $\omega\subset\R^2$ with a sufficiently
regular boundary $\partial\omega$ be given. The problem consists in
finding a set $\Omega\supset\bomega$ 
and a function $u:\Omega\setminus\bomega\to\R$ satisfying
\begin{equation}
\left\{\begin{aligned}
  \Delta u =0 &\quad\text{in }\Omega\setminus\bomega\\
  u=1 &\quad\text{on } \partial\omega \\
  \left.\begin{aligned}&u=0 \\
        &\PD u\nvec=\gamma
        \end{aligned}\right\} &\quad\text{on }\partial\Omega.
\end{aligned}\right.
\tag{$\widetilde\calP(\omega)$}
\label{eq:bernoulli}
\end{equation}

This paper deals with the control of the shape of $\Omega$
in \eqref{eq:bernoulli}. The geometry of the free boundary $\partial \Omega$
will be driven by the shape of $\omega$ towards a given target free
boundary $\Gamma_t$. 
Our optimization problem reads as follows:
\begin{equation}
\left\{\begin{aligned}
 \null &\text{Find }\omega^*\in\calOt \text{ such that}\\
       &J(\Gamma(\omega^*)) \le J(\Gamma(\omega))
       \end{aligned}\right.
 \tag{$\widetilde{\mathbb P}$}
\end{equation}
holds for any $\omega\in\calOt$, where $\calOt$ is a set of admissible designs.
The cost functional to be considered is the distance between
$\Gamma_t$ and the free boundary $\Gamma(\omega) := \partial
\Omega(\omega)$ 
being the solution of \eqref{eq:bernoulli}, i.e.
\begin{equation}
J(\Gamma(\omega)) = \rho(\Gamma(\omega), \Gamma_t),
 \label{eq:objf}
\end{equation}
where $\rho$ is a function characterizing the distance of $\Gamma(\omega)$
from $\Gamma_t$.

In this paper we make the assumption that $\Omega$ is a {\em star-like domain}.
However, in contrast to \cite{toi08} no such assumption is made on
$\omega$ meaning that also the topology of $\omega$ may change.

\section{Pseudo-solid formulation of Bernoulli free boundary problems}
\label{sec:pseudosolid}

Free boundary problems have in common the difficulty that the geometry
(here the domain $\Omega$) has to be determined simultaneously with
the solution $u$ of the state problem, which implies that a numerical
solution has to be done iteratively \cite{kar99}.  Possible solution
strategies include trial methods, linearization methods (continuous or
discrete) \cite{cuv90}, and shape optimization methods
\cite{has03}. We decided for  the
so-called {\em pseudo-solid approach} (PSA) in which the unknown
domain $\Omega$ is obtained by an appropriate deformation of the
reference configuration $\hatOmega$. This deformation is one of
the unknowns in PSA.
The main advantage of this approach
is that there is no need to construct an explicit parametrization
of the geometry using e.g. a conformal mapping of the reference domain.
PSA is useful also in numerical
realization: finite element partitions of $\Omega$  can be constructed
via the respective deformation of the partition of $\hatOmega$.
As the Bernoulli problem has to be solved  several times for different
$\omega$, the choice of an efficient solver is also important. We use
Newton's method because of its fast convergence.
Moreover, this solution strategy readily allows us to
obtain geometrical sensitivities of the system.

Let $\hatOmega\subset\R^2$ be a fixed, simply connected reference
domain. In the
pseudo-solid technique we construct a mapping $F:\R^2\to\R^2, \
\hatOmega\mapsto F(\hatOmega)=:\Omega$ such that $\Omega$ solves
\eqref{eq:bernoulli} for given $\omega$.  To construct such
$F$ we treat $\hatOmega$ as an elastic solid that undergoes a
deformation caused by an external loading $p$ such that the deformed
solid defines such 
$\Omega$. Thus, problem \eqref{eq:bernoulli} is strongly coupled with
the Lam\'e system of linear elasticity in which the loading $p$ applied to
$\partial\hatOmega$ 
plays the role of an unknown in PSA.
This approach has been previously used to solve free surface
flow problems (see e.g.\ \cite{cai00, sou01}) and Bernoulli free boundary
problems in \cite{toi08}.

For any $\ww\in
\Wad:=$\{''sufficiently'' small and regular deformations\} we define the
domain
$$
  \Omega_\ww = \{ \xx\in\R^2 \mid \xx=\hat\xx + \ww(\hat\xx),
  \quad\hat\xx\in\hatOmega \}.
$$
We introduce the following function spaces:
\begin{align*}
W_\omega &=\{ \ww\in[H^1(\hatOmega)]^2 \mid
\ww|_{\omega}=\0 \} \\
V_c(\Omega) &=\{ \varphi\in H^1(\Omega) \mid \varphi|_\omega = c \},
\quad c\in\R.
\end{align*}
The pseudo-solid formulation of our free boundary
problem then reads as follows: Given $\omega$,\\
find
$(u,p,\vv)\in V_1(\Omega_\vv)\times L^2(\partial \hatOmega) \times W_\omega$ such
that
\begin{equation}
\left\{
\begin{aligned}
  \null &\int_{\Omega_\vv \setminus \omega} \nabla
  u\cdot\nabla\varphi\,d\xx = \gamma\int_{\partial \Omega_\vv} \varphi\,ds \quad  \forall
  \varphi\in V_0(\Omega_\vv) \\
  \null &\int_{\partial \Omega_\vv} u\psi\,ds = 0 \quad \forall \psi\in L^2(\partial \Omega_\vv)
  \\
  \null &\int_{\hatOmega \setminus \omega}
  \sigma(\vv):\varepsilon(\ww)\,d\xx = \int_{\partial \hatOmega}
  p\,\nvec\cdot\ww\,ds \quad \forall \ww\in W_\omega.
\end{aligned}
\right.
\tag{$\calP(\omega)$}
\label{eq:bernoulli-1}
\end{equation}

Equations $(\calP(\omega)_1)$ and $(\calP(\omega)_2)$ constitute the
weak form of \eqref{eq:bernoulli} while
$(\calP(\omega)_3)$ is the weak form of the linear elasticity
problem in $\hatOmega \setminus \bomega$. Here $p\nvec$ is an
(unknown) external load. The 
components of the strain and stress tensors
$\varepsilon=\{\varepsilon_{ij}(\vv)\}$and
$\sigma=\{\sigma_{ij}(\vv)\}$ associated with a displacement field
$\vv$ are given by
$$
\varepsilon_{ij}(\vv)=\frac12\left( \frac{\partial v_i}{\partial
x_j} + \frac{\partial v_j}{\partial x_i} \right), \quad
  \sigma_{ij}(\vv)=2\mu\epsilon_{ij}(\vv)+ \lambda \delta_{ij} \nabla \cdot \vv,
\quad i,j=1,2,
$$
respectively, where $\mu$ and $\lambda$ are
the Lam\'e coefficients. Since in this case the linear elasticity
system does not have any physical meaning, the Lam\'e coefficients
can be chosen quite freely. In this paper the choice $\mu = 0.5$
and $\lambda = 0$ was made. The solvability of \eqref{eq:bernoulli-1} is
analyzed in details in \cite{toi08}. The relation between
\eqref{eq:bernoulli} and \eqref{eq:bernoulli-1} is readily seen: if
$(u,p,\vv)$ is a solution of \eqref{eq:bernoulli-1} then the couple
$(u_{\big|\Omega_\vv\setminus\bomega},\,\Omega_\vv)$ solves \eqref{eq:bernoulli}.

 \begin{figure}[h]
 \center
 \includegraphics[width=0.9\textwidth]{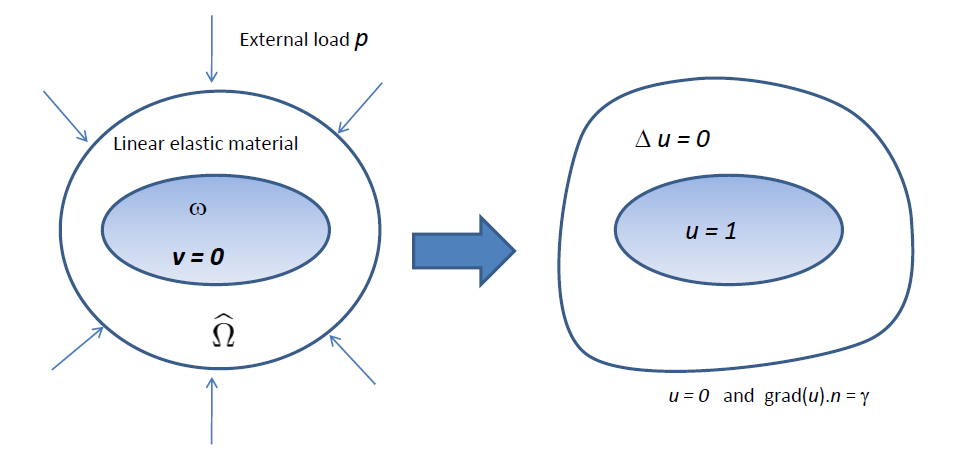}

\null\hfill a) reference configuration \hfill b) deformed configuration\hfill\null
 \caption{Principle of the pseudo-solid approach.}
 \label{fig:principle}
 \end{figure}

\section{Level set approach to the optimization problem}
In the previous paper \cite{toi08} the design domain $\omega$ was
parameterized using the polar co-ordinates. A tendency towards
fractal-like designs indicating possible 
topological changes of $\omega$ was observed in certain cases.
However, the boundary variation technique used in \cite{toi08} is not
able to handle topological changes automatically. This will be done by
 a level set parametrization \cite{osh88,osh03} of $\omega$.

The basic idea is simple:
Let $D$ be a larger domain containing all admissible $\omega$ (for example
a rectangle slightly larger than the
bounding box of the target boundary $\Gamma_t$).
Let $\psi:D\to\R$ be given and define the set $\omega$ by
\begin{equation}
\omega := \omega(\psi) :=\left \{ \xx \in D \, | \, \psi(\xx) > 0
\right \}, \quad 
\psi \in \Uad.
 \label{eq:omega-psi-def}
\end{equation}
Here $\Uad$ is a family of admissible level set functions 
such that $\psi\in\Uad \implies
\omega(\psi)\in\tilde\calO$. Clearly the parametrization
\eqref{eq:omega-psi-def} allows topological changes of $\omega$.

We can now reformulate problem $(\widetilde\PP)$ as follows:
\begin{equation}
\left\{\begin{aligned}
 \null &\text{Find }\psi^*\in\Uad \text{ such that}\\
       &J(\Gamma(\psi^*)) \le J(\Gamma(\psi)) \quad \forall\psi\in\Uad,
       \end{aligned}\right.
 \tag{$\widehat{\PP}$}
\end{equation}
where $\Gamma(\psi)$ is the free boundary defined by $(\calP(\omega(\psi)))$

Next we introduce the following relaxed state problem which does not
contain explicitly any Dirichlet type boundary conditions. For $\eps>0, \
\psi\in\Uad$ given:\\ 
Find
$(u_\eps,p_\eps,\vv_\eps)\in H^1(\Omega_{\vv_\eps})\times L^2(\partial\hatOmega)\times [H^1(\hatOmega)]^2$
such that
\begin{equation}
\left\{\begin{aligned}
  \null &\int_{\Omega_{\vv_\eps}} \nabla
  u_\eps\cdot\nabla\varphi\,d\xx - \gamma\int_{\partial \Omega_{\vv_\eps}}
  \varphi\,ds + \int_{\Omega_{\vv_\eps}}
  G_\epsilon(\psi)(u_\eps-1)\varphi\,d\xx = 0 \quad  \forall 
  \varphi\in H^1(\Omega_{\vv_\eps}) \\
  \null &\int_{\partial \Omega_{\vv_\eps}} u_\eps\psi\,ds = 0 \quad
  \forall \psi\in L^2(\partial \Omega_{\vv_\eps}) 
  \\
  \null &\int_{\hatOmega}
  \sigma(\vv_\eps):\varepsilon(\ww)\,d\xx - \int_{\partial \hatOmega}
  p_\eps\,\nvec\cdot\ww\,ds + \int_{\hatOmega} G_\epsilon(\psi)
  \vv_\eps \cdot \ww\,d\xx = 0 
  \quad \forall \ww\in [H^1(\hatOmega)]^2.
\end{aligned}\right.
 \tag{$\calP_\epsilon(\psi)$}
 \label{eq:bernoulli-p-e}
\end{equation}
Above, $G_\eps$ is a penalty function releasing the constraints $u=0$
and $\vv=\0$ in $\omega$. The classical choice is    
$G_\eps=\frac1\epsilon H$, where $H$ is the 
Heaviside function. This choice of $G_\eps$ will be used in what follows.

We define now the ``relaxed'' optimization problem
\begin{equation}
\left\{\begin{aligned}
\null&\text{Find } \psi_\epsilon^\star \in\Uad \text{ such that}\\
&J(\Gamma_\eps(\psi^\star_\epsilon) \le J(\Gamma_\eps(\psi)) \quad
\forall \psi\in\Uad,
\end{aligned}\right.
\tag{$\PP_\epsilon$}
\end{equation}
where $\Gamma_\eps(\psi)=F_\eps(\partial\hatOmega), \
F_\eps=\text{id}+\vv_\eps$ and $\vv_\eps$ is the third component of
the solution of \eqref{eq:bernoulli-p-e}.


\section{Discretization}
\label{sec:discr} 
One of the advantages of level set methods is that they
avoid tracking of the boundary of the design domain. Instead a
fixed mesh is used. In our case the mesh $\hat\calT$ of the reference domain
$\hatOmega$ is fixed, but the mesh $\calT$ of $\Omega$ is moving according
to the pseudo-solid strategy. However, the boundary of $\omega$ is
not exactly tracked in either of these meshes. An unstructured mesh consisting
of triangles is used to approximate $\hatOmega$. The mesh of
$\Omega_\vv$ is obtained by displacing the nodes of
$\hat\calT$ using the discrete displacement field
$\vv_h$ which approximates $\vv$.

\subsection{Discretization of the state problem}

Let $\eps>0,\ \psi\in\Uad$ be given. To simplify notation, the penalty
parameter $\eps$ at the discrete solution will be omitted.
The finite element discretization of \eqref{eq:bernoulli-p-e} reads as
follows:\\
Find
$(u_h,p_h,\vv_h)\in V_h(\Omega_{\vv_h})\times Q_h\times \WW_h$
such that
\begin{equation}
\left\{
\begin{aligned}
  \null &\int_{\Omega_{\vv_h}}\!\! \nabla
  u_h\cdot\nabla\varphi_h\,d\xx -
\gamma\int_{\partial \Omega_{\vv_h}} \!\!\!\!\!\!\varphi_h\,ds
      + \int_{\Omega_{\vv_h}}
      \!\!\!\!\!G_\eps(\psi)(u_h-1)\varphi_h\,d\xx = 0 \quad  \forall 
  \varphi_h\in V_h(\Omega_{\vv_h}) \\
  \null &\int_{\partial \Omega_{\vv_h}} u_h\mu_h\,ds = 0
              \quad \forall \mu_h\in Q_h(\partial\Omega_{\vv_h}) :=
      \left.V_h(\Omega_{\vv_h})\right|_{\partial\Omega_{\vv_h}}
  \\
  \null &\int_{\hatOmega}
  \sigma(\vv_h):\varepsilon(\ww_h)\,d\xx - \int_{\partial \hatOmega}
  p_h\,\nvec\cdot\ww_h\,ds + \int_{\hatOmega} G_\eps(\psi) \vv_h \cdot
  \ww_h\,d\xx = 0 
            \quad \forall \ww_h\in W_h,
\end{aligned}
\right.
 \tag{$\calP_\eps^h(\psi)$}
\label{eq:bernoulli-p-e-h}
\end{equation}
where $V_h(\Omega_{\vv_h}), \ \WW_h=W_h\times W_h,\ Q_h$ are finite element 
approximations of $H^1(\Omega_\vv)$, 
$[H^1(\hatOmega)]^2$ and $L^2(\partial\hatOmega)$ respectively. Here we shall use linear
triangular elements for constructing $V_h(\Omega_{\vv_h})$ and $W_h$, while 
$Q_h=\left.W_h\right|_{\partial\hatOmega}$.

\subsection{Discrete optimization problem}
In the traditional level set method the function $\psi$
is taken to be a function of pseudo-time $t$, $\psi := \psi(\xx,t)$,
and the optimization process is realized by solving the
Hamilton-Jacobi equation 
\begin{equation}
\frac{\partial \psi}{\partial t} + v_n \| \nabla \psi \| = 0.
\end{equation}
Here $v_n$ is the velocity derived by the
means of sensitivity analysis, often done on the continuous level.
The function $\psi$ is then advanced towards the steady state in 
pseudo-time, see e.g.\ \cite{all04}.

Despite of the conceptual simplicity it is not so straightforward to
implement the conventional level set 
method  due to the need of
appropriate upwind schemes, an extension of the velocities and
re-initialization algorithms. Indeed, since the Hamilton-Jacobi equation
does not in general admit a smooth solution, an appropriate
upwind scheme must be used for the time integration. The velocity
$v_n$ is often meaningful only on the boundary $\partial \omega$,
and must be extended to the whole domain, or at least into a
neighborhood of $\partial \omega$. Finally, the function
$\psi$ should be an approximation of the signed distance function, i.e.
$\psi(\xx) \approx \text{sign}(\psi(\xx)) || \xx - \xx_0 ||$, where
$\xx_0$ is the 
closest point to $\xx$ for which $\psi(\xx_0) = 0$. To force this
property a re-initialization procedure is often used.

Several approaches to overcome these difficulties
have been proposed.
In \cite{wan06} radial basis functions (RBF) are used to
define the function $\psi$, and the Hamilton-Jacobi
equation is transformed into a system of ordinary differential
equations. In \cite{che07} the function $\psi$ is
constructed combining parameterized primitives with a
radial basis function representation of the so called
R-functions. In \cite{bel03}, \cite{nei09} the function $\psi$ is
approximated by the same shape functions on the same mesh used to solve the
state problem.
In this paper we follow \cite{luo08} and utilize the compactly supported 
$C^2$-continuous radial basis  functions \cite{wen05} to 
parameterize explicitly the level set function. Then the geometry
will inherit  
smoothness properties of the underlying parameterization.
Moreover, we have a fixed set of design variables, and one can
use sophisticated optimization methods instead of 
performing integration in the pseudo-time.

\begin{figure}[h]
\center
\includegraphics[width=0.9\textwidth]{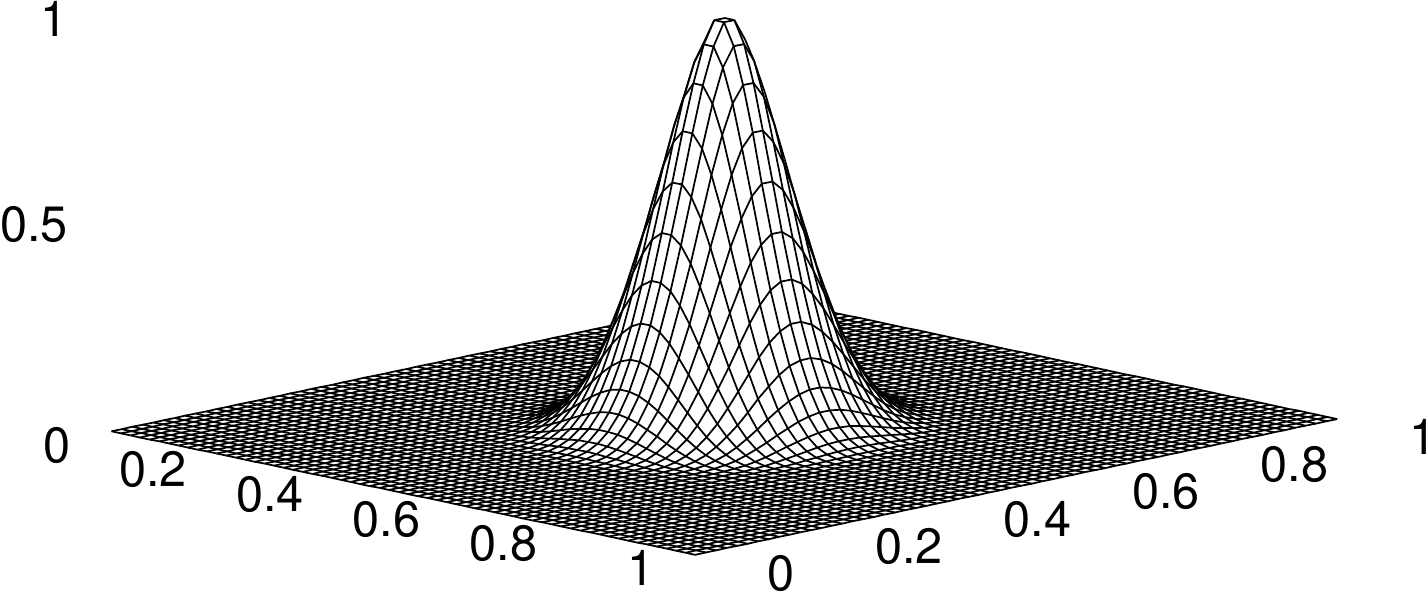}
\bigskip
\caption{Surface plot of a radial basis function $\psi_{ij}$}
 \label{fig:RBF}
\end{figure}

We introduce a set of $N \times N$ basis functions, whose knots are placed
in the interior of the domain $D$ as follows.
The coordinates of the knot $(i,j)$ are given by
\begin{align}
x_{ij} &= x_{min} + (j-1) \frac{x_{max}-x_{min}}{N-1}, \quad i,j=1,...,N \\
y_{ij} &= y_{min} + (i-1) \frac{y_{max}-y_{min}}{N-1}, \quad i,j=1,...,N
\end{align}
where $x_{min}$ and $x_{max}$ are the minimal and maximal
$x$-coordinates of the rectangle $D$, respectively and similarly
$y_{min}$ and $y_{max}$. 
The  RBF associated with this knot is then
\begin{equation}
  \psi_{ij}(r_{ij}) = \max\left\{0, (1-r_{ij})\right\}^4 (4r_{ij}+1),
\end{equation}
where
\begin{equation}
r_{ij} = \frac{\sqrt{ (x-x_{ij})^2 + (y-y_{ij})^2  }}{r_s}.
\end{equation}
Here $r_s>0$ is a given parameter, the radius of the support (see
Figure~\ref{fig:RBF}). 
To guarantee the overlapping of the supports of $\psi_{ij}$ we define this
parameter as
\begin{equation}
r_s = 4\cdot\max \left( \frac{ x_{max}-x_{min} }{N-1}, \frac{ y_{max}-y_{min} }{N-1} \right).
\end{equation}

The level set function $\psi$ is then approximated by the linear
combination of ${\psi_{ij}}$:
\begin{equation}
\psi:=\psi_N(\Alpha) = \sum_{i,j=1}^N \alpha_{ij} \psi_{ij}.
 \label{eq:linear-comb}
\end{equation}
Thus, the discrete design variables of the parameterized optimization
problem 
are represented by the vector
$\bld \alpha = (\alpha_{11}, \alpha_{12}, \ldots ,\alpha_{NN})$. 

Using the assumption that $\Omega$ is a star-like domain, 
the objective functional \eqref{eq:objf} will be
 given in the discrete setting by
\begin{equation}
 \calJ(\Alpha) := J(\vv_h(\Alpha))
=\int_{0}^{2 \pi} (g_\Alpha(\theta) - g_t(\theta))^2 \, d\theta,
\label{eq:discr-objf}
\end{equation}
where $\Alpha$ is the vector of the discrete design variables,
$\vv_h(\Alpha)$ is a part of the solution to
$(\calP^h_\eps(\psi_N(\Alpha)))$, 
$g_\Alpha(\theta)$ is the radius of the free boundary
corresponding to $\Alpha$, and $g_t(\theta)$ is the radius of
the target boundary at the angle $\theta$. The set $\Uad$ of all
admissible level  set functions is represented by functions of the
form \eqref{eq:linear-comb} with
$\Alpha\in\calU:=[\alpha_{\text{min}},\alpha_{\text{max}}]^{N\times N}
\subset\R^{N\times N}$, with $\alpha_{\text{min}},\alpha_{\text{max}}$ given.

Thus, the finite-dimensional optimization problem to be realized
reads as follows:

\begin{equation}
\left\{
\begin{aligned}
 \null &\text{Find } \Alpha^*\in\calU \quad\text{such that}\\
 &\calJ(\Alpha^*) \le \calJ(\Alpha) \quad \forall\Alpha\in\calU
\end{aligned}
\right.
\end{equation}
subject to $(\calP^h_\eps(\psi_N(\Alpha)))$. 

\subsection{Construction of the initial guess}
Obviously, the location of the free boundary corresponding to a given inner
boundary $\partial\omega$ is not known a priori. Therefore we start
the optimization from 
a simple configuration where the location of the free boundary
corresponding to the initial $\omega$ is known.

Our initial guess of $\Alpha$ is constructed in such a way that 
\begin{equation}
\psi_N(x_{ij}, y_{ij}) = \hat{\psi}_1 (x_{ij}, y_{ij}) \qquad \forall \,
i,j = 1,\ldots,N
\label{eq:init-guess-system}
\end{equation}
where $(x_{ij}, y_{ij})$ are the knots of the radial basis
functions and $\hat{\psi}_1$ is given by \eqref{eq:exscalar} in
Appendix~2 with
$R=1$. We obtain a system of linear equations, which is
known to be invertible. The initial values for the
design vector $\Alpha$ are solutions to
\eqref{eq:init-guess-system}. 
The zero level set of such function $\psi$ then
approximates $\partial B(0,1)$.

For given $\gamma$, the initial reference domain
$\hatOmega$ is taken to be the circle $B(0, C(1, \gamma))$ (for the
definition of $C(1,\gamma)$ see Appendix~2).

\subsection{Algebraic form and Sensitivity analysis}

Let $\psi\in U_{ad}$ be given.
In the pseudo-solid approach we simultaneously seek the scalar
function $u$, the pressure $p$, and the deformation field $\vv$ 
which deforms the
reference domain $\hatOmega$ into the one that solves
\eqref{eq:bernoulli}. The elasticity system
\eqref{eq:bernoulli-p-e-h}${}_3$ is 
thus solved in the undeformed configuration $\hatOmega$ 
of the pseudo-solid, whereas
equations \eqref{eq:bernoulli-p-e-h}${}_1$ and
\eqref{eq:bernoulli-p-e-h}${}_2$ are 
solved in the deformed one. 
Therefore, they have to be discretized 
by different meshes, too. Let us denote the nodal co-ordinates of
the triangulation $\hat\calT$ of $\hatOmega$ by $\hat\XX$. We
simply transform this triangulation into the one of $\Omega_{\vv_h}$
which is characterized by the nodes 
$\XX=F_h(\hat\XX)$, where $F_h$ is defined by the discrete
displacement field $\vv_h$ being the approximation of $\vv$, i.e.
\begin{equation}
\label{eq:mesh_relation} \XX = \hat\XX + \vv_h(\hat\XX).
\end{equation}

The algebraic form  of 
\eqref{eq:bernoulli-p-e-h}${}_1$, \eqref{eq:bernoulli-p-e-h}${}_2$
resulting from 
an appropriate discretization can be written  as 
\mbox{$\rr_1(\qq_u,\qq_v)=\0$}, $\rr_2(\qq_u,\qq_v)=\0$, respectively
 and the linear elasticity system 
\eqref{eq:bernoulli-p-e-h}${}_3$  as $\rr_3(\qq_v, \qq_p) = \0$, where
$\qq_u$, $\qq_v$ and $\qq_p$ are the nodal values of $u_h$, $\vv_h$ and
$p_h$, respectively. Here the
dependence of $\rr_1$ and $\rr_2$ on $\qq_v$ is through the 
nodal co-ordinates, as specified by 
\eqref{eq:mesh_relation}. Dimensions of the vectors $\qq_u$,
$\qq_v$ and $\qq_p$ are $n$, $2n$ and $n_e$ respectively, where
$n$ is the number of the nodes in $\hat\calT$ and $n_e$ is the number of
the nodes on $\partial\hatOmega$.

Let us introduce the following notation:
\begin{equation}
\label{eq:q_r} \qq =
\begin{bmatrix}
\qq_u\\
\qq_v\\
\qq_p\\
\end{bmatrix}
 \quad \text{ and } \quad \rr =
\begin{bmatrix}
\rr_1\\
\rr_2\\
\rr_3
\end{bmatrix}.
\end{equation}
Then the algebraic form of the discretized coupled system
\eqref{eq:bernoulli-p-e-h} can be written in short
as $\rr(\qq)=\0$. This system will be solved using Newton's method:
\begin{equation}
\label{eq:newton_short} \qq^{(k+1)}=\qq^{(k)}-\left(
\frac{\partial \rr(\qq^{(k)})}{\partial \qq}\right)^{-1}
\rr(\qq^{(k)})
\end{equation}
with the Jacobian matrix 
\begin{equation}
\label{eq:jacobian} \left( \frac{\partial \rr}{\partial \qq}
\right) =
\begin{bmatrix}
\frac{\partial \rr_1}{\partial \qq_u} & \frac{\partial \rr_1}{\partial \qq_v} & \0 \\
\noalign{\medskip}
\frac{\partial \rr_2}{\partial \qq_u} & \frac{\partial \rr_2}{\partial \qq_v} & \0 \\
\noalign{\medskip} \0 & \frac{\partial \rr_3}{\partial \qq_v}
&\frac{\partial \rr_3}{\partial \qq_p}
\end{bmatrix}
\end{equation}
%
%
%


\begin{remark} \rm
The discrete cost function $\calJ$ defined by \eqref{eq:discr-objf} is
not differentiable due to the 
discontinuous Heaviside function $H$. If one wishes to use descent type
optimization methods (and this will be our case),  then
smoothing of the Heaviside function is necessary (see Appendix~1) to make
the components 
of the residual vector $\rr$ continuously differentiable functions of the
design variables. Thus, while assembling the discrete system arising
from \eqref{eq:bernoulli-p-e-h}, 
$H$ will be replaced by a smoothed function $H_\beta$,
i.e.\ $G_\eps=\frac1\eps H_\beta$ in \eqref{eq:bernoulli-p-e-h} in
what follows.
\end{remark}

\begin{remark} \rm
Notice that the equations are coupled in a quite complicated way.
For example, the residual $\rr_1$ depends naturally on $\qq_u$, but
also on $\qq_v$ through the shapes
of the elements as specified by \eqref{eq:mesh_relation}.
Moreover, since the mesh is moving and $\psi$ is a function
of location, $\rr_1$ depends on $\qq_v$
also through $H_\beta(\psi)$. This dependence must be taken
into account especially during the sensitivity analysis phase
in order to obtain perfectly consistent derivatives.

Despite this nonstandard coupling between the equations,
the Jacobian matrix \eqref{eq:jacobian} is easy to compute using
the sparse forward mode automatic differentiation \cite{bis96-2}.
Our implementation of the automatic differentiation technique is described in 
\cite{toi10-2}.
For a general introduction to the principles of
automatic differentiation see \cite{gri08}.
\end{remark}

\def\kmax{{k_{\text{max}}}}
\def\imax{{i_{\text{max}}}}

If $\calJ$ is smooth then using the well-known adjoint approach, the gradient
$\nabla_\Alpha \calJ$ can be computed from
\begin{equation}
\nabla_\Alpha \calJ = \left(\frac{\partial \rr}{\partial
\Alpha}\right)^T \boldsymbol\nu,
 \label{eq:adj-1}
\end{equation}
where the adjoint vector $\boldsymbol\nu$ solves the adjoint equation
\begin{equation}
\label{eq:adj-2}
  \left(\frac{\partial \rr}{\partial \qq}\right)^T \boldsymbol\nu =
  - \nabla_\qq\calJ.
\end{equation}
The required Jacobian matrices and gradients in \eqref{eq:adj-1} and
\eqref{eq:adj-2} can be again easily computed
using the tools of automatic differentiation.
Notice that the Jacobian $\partial \rr / \partial \Alpha$
has a sparse structure, since the radial basis functions are compactly
supported. This sparsity is automatically exploited, since
we use the sparse forward mode automatic differentiation
technique. To avoid going through all radial basis
functions while evaluating $\psi$, a quadtree data structure is exploited
to exclude RBFs that can not have a non-zero value at
the point of evaluation.

\subsection{Optimization strategy using remeshing}
As explained in the previous paper \cite{toi08}, remeshing is sometimes
needed since the mesh deformation approach can not handle too large
displacements. Indeed, the Newton method used to solve the coupled system 
may not converge in case of excessive deformations of $\hat\Omega$. Moreover,
if the mesh gets too distorted, significant errors in the solution and
numerical instabilities may appear.

In this paper we adopt the following adaptive optimization strategy: 
to avoid too large
deformation fields $\vv_h$ in \eqref{eq:bernoulli-p-e-h} the reference
domain $\hatOmega$ is re-initialized 
after $\kmax$ optimization steps. 
Moreover, the half width $\delta$ of the gray region related to the
smoothing of the Heaviside 
function (see Appendix~1) is determined adaptively. 
This $\delta$ is then used to determine the smoothing parameter
$\beta$, i.e.\ $\beta:=\beta(\delta)$. Note that due to this adaptive
strategy $\beta$ varies 
spatially, too. Therefore there are two smoothing parameters $\eps,\delta$
related to the discrete pseudo-solid problem which we thus denote by
$(\calP_{\eps,\delta}^h(\psi))$.

The strategy of
choosing $\delta$  is also altered after $\kmax$  optimization steps. In early
steps a larger value  of $\delta$
is used, which makes the objective function smoother
and enables fast progress. The value of $\delta$ is then gradually decreased,
improving the approximation of the exact Heaviside step function.
We choose the half width $\delta$ of the desired gray region
after $i$:th re-initialization step to be
\begin{equation}
\delta = \frac12(i_{\text{max}}-i+1) h,
\label{eq:delta}
\end{equation}
where $i_{\text{max}}$ is the number of re-initializations to be done
and $h$ is the characteristic mesh size.

We choose the number of re-initializations to be modest, e.g.\  $\imax=8$.
In the early stages of optimization the progress is rapid, and changes
in the design domain $\omega$ are large. 
Therefore we let the number of steps between re-initializations to
increase in the course of optimization, e.g.\
$k_{\text{max}}=5 \cdot 2^i$.

The re-initialized
optimization process can be described as the following algorithm:

\begin{align}
&\text{{\bf Algorithm 1}}\quad (\imax,h,toler \ \text{ given}) \nonumber \\
&\text{Compute initial guess }\Alpha \text{ using
    eq.\ \eqref{eq:init-guess-system}} \ ; \\ 
&\text{Set } \hatOmega := B(0,C(1,\gamma)) \ ;   \label{alg:eka}\\
&\text{{\bf for }}i=1,...,i_{\text{max}} \\
&\quad \text{Generate triangulation }\hat\calT \text{ of }\hatOmega \ ;\\
&\quad \text{Set }\delta:=\frac12(i_{\text{max}}-i+1)h \text{ and }
  \kmax=5\cdot 2^i\ ;\\
&\quad \text{{\bf do }} \\
&\qquad \text{Set }k:=1 \ ;\\
&\qquad \text{Solve state problem } (\calP^h_{\eps,\delta}(\psi_N(\Alpha)))
  \ ;\\
&\qquad \text{Evaluate }
  \cal J(\Alpha) \text{ and }\nabla_\Alpha \calJ(\Alpha) 
  \text{ using }\eqref{eq:discr-objf},\ \eqref{eq:adj-1}
          \text{ and  }\eqref{eq:adj-2} ;\\   
&\qquad \text{Find descent direction } \tilde\Alpha \text{ of
          }\calJ \text{ at } \Alpha \ ; \label{eq:opti-step}\\
&\qquad \text{Update design } \Alpha := \Alpha + \tilde\Alpha \ ;\\
&\qquad \text{Set } k:=k+1 \ ;\\
&\quad \text{{\bf while }} k<k_{\text{max}} \text{ and }
  \|\tilde\Alpha\|>toler \ ;\\
&\quad \text{Reinitialize }\hatOmega:=\Omega_{\vv_h} \text{ using
    the latest } \vv_h \ ;\label{eq:re-ini}\\
&\text{{\bf end}} \label{alg:vika}
\end{align}


The re-initialization of $\hatOmega$ in step \eqref{eq:re-ini} is done
as follows.  
One fits by least squares a cubic B-spline curve 
to the outer boundary of the deformed mesh $\calT(\vv_h)$
obtained as the solution to the free boundary problem
corresponding to the current design $\Alpha$. 
Nodes are then distributed on this curve, and this geometry is 
given to the mesh generator Gmsh \cite{geu09}.

Thus, in algorithm \eqref{alg:eka}--\eqref{alg:vika} we
basically solve $\imax$ different optimization problems, 
performing a mesh regeneration for each problem.
The initial guess for the next problem is always the best design found
so far, i.e.\ the end result of the previous problem.

\section{Numerical examples} \label{sec:examples}

In this section we illustrate the performance of the proposed method.
The target domains are the
same as in \cite{toi08}. The value of the penalty parameter $\epsilon$ 
is $10^{-3}$. 
We used the gradient based optimizer Donlp2 \cite{donlp2} to
realize the step \eqref{eq:opti-step} in Algorithm~1. 
The parameters
  $\alpha_{\text{min}},\alpha_{\text{max}}$ were chosen as $\pm10^{20}$,
i.e.\ we have practically unconstrained minimization problem. The
performed numerical computations indicated no need
to pose more strict constraints.

\begin{example}\label{ex:rsq-1} \rm

The target $\Gamma_t$ is the ``rounded square''. The length of sides of
the square is 4. Each corner is rounded using a quarter of a circle
of radius 1. For
the magnitude of the normal derivative the value $\gamma=-1$ was
used. It was observed in \cite{toi08} that if the family $\calOt$ of
admissible inner inclusions contains only star-like domains $\omega$,
then the target $\Gamma_t$ can never be matched. Moreover, if the number
of the design variables increased,
the boundary $\partial\omega$ became more and more oscillating.

The final reference domain $\hatOmega$ was discretized using 18842 linear
triangular elements. 
The number of RBFs used was $30\times30$. The initial value
of the cost functional was 1.09. After  630 optimization steps (and
2376 function evaluations) the cost function reduced to the value $3.78
\cdot 10^{-5}$. 

The final  geometry and the contour plot of $H_\beta(\psi_N^\star)$ are shown in
Figure~\ref{fig:rsq-1}. 
\begin{figure}[h]
\center
\includegraphics[width=0.45\textwidth]{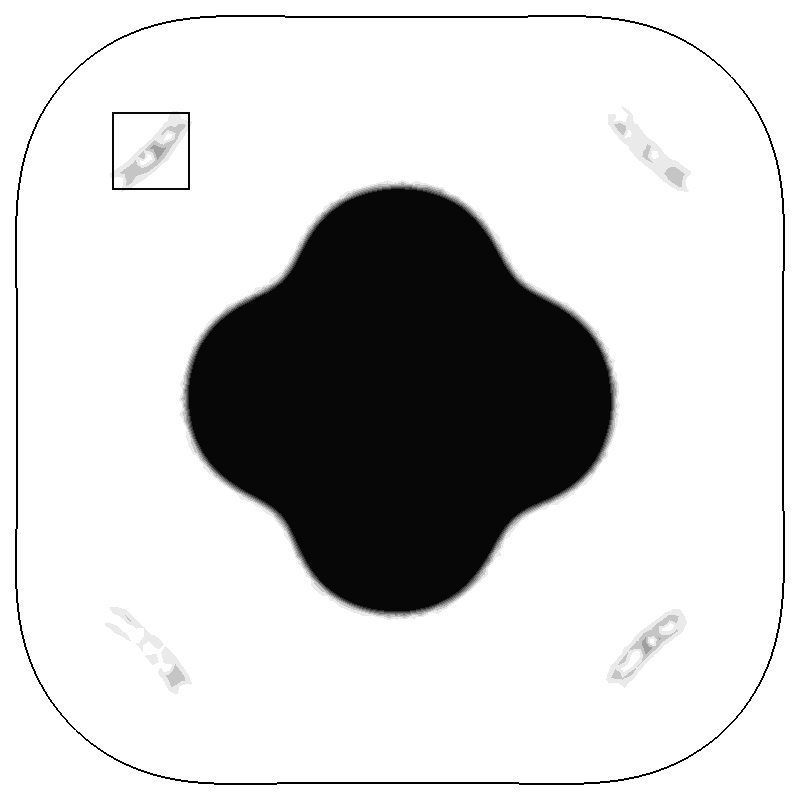}
\includegraphics[width=0.45\textwidth]{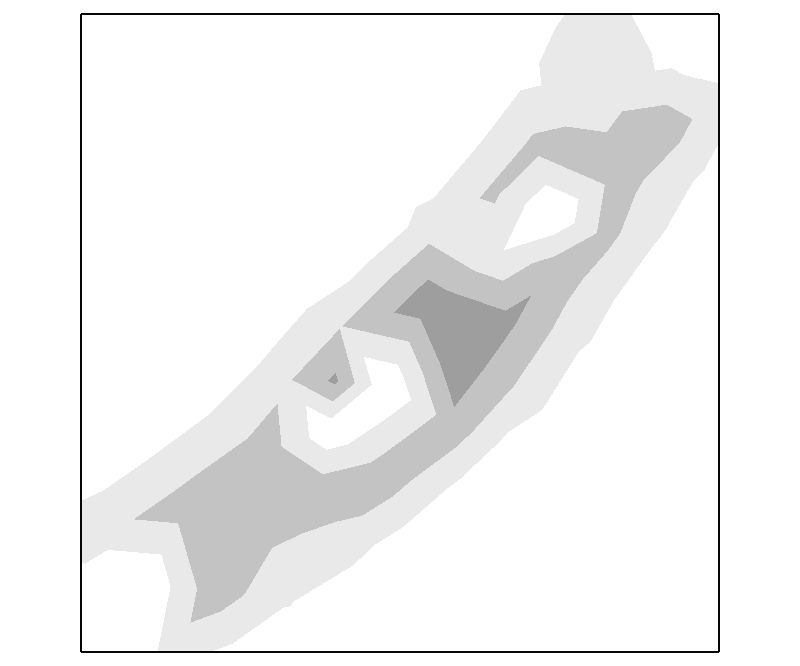}
\caption{Left: the final geometry and $H_\beta(\psi_N^\star)$. Right:
  zoom of the framed area.}
 \label{fig:rsq-1}
\end{figure}

In addition to ``black'' (i.e.\ $\omega$) and ``white''
(i.e.\ $\Omega\setminus\overline{\omega}$) regions,
 smoothing of the Heaviside function also
produces by its construction ``grey'' regions (where
$0<H_\beta(\psi_N^\star)<1$) near the 
zero level set of $\psi$ representing 
$\partial\omega$. In the final design however grey 
regions which can not be interpreted as being close to
$\partial\omega$ may appear. 
In Figure~\ref{fig:rsq-1} we see one of those problematic 
regions. By examining  the potential $u_h$ (see Figure~\ref{fig:rsq-1-pot})
in that area,
we find that the maximum value of $u_h$ in this region is
only about $0.86$, whereas in $\omega$ we should approximately meet
the condition $u_h=1$.

\begin{figure}[h]
{\center
\includegraphics[width=0.4\textwidth]{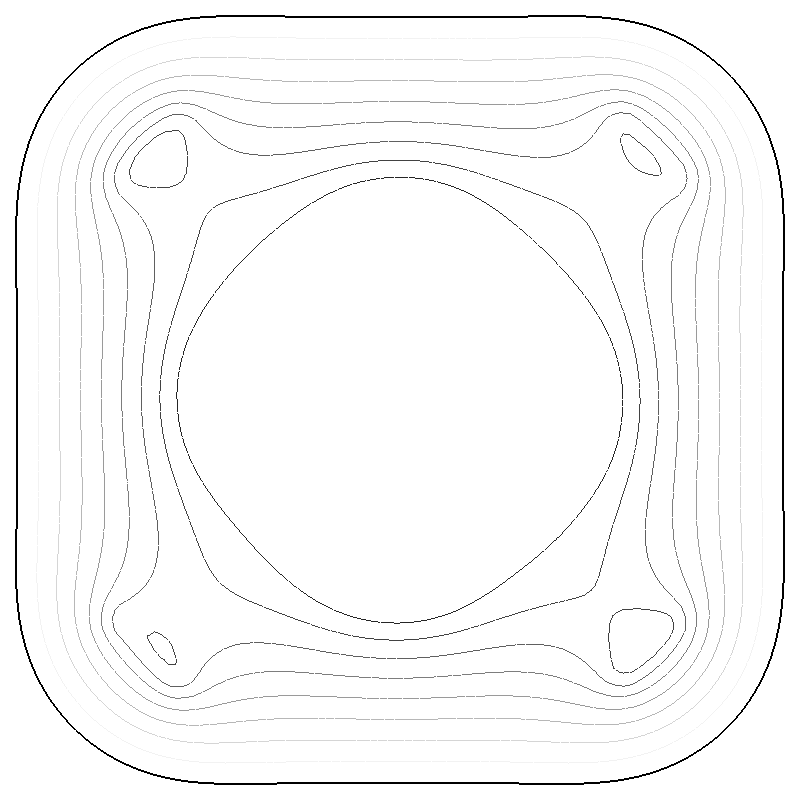}\\
}
\caption{Contour plot of the potential $u_h$ related to Example \ref{ex:rsq-1}.}
\label{fig:rsq-1-pot}
\end{figure}
\end{example}

\begin{example}\label{ex:rsq-2} \rm
To prevent such problematic grey regions we propose to add the
following penalty term to the cost functional $\calJ$:
\begin{equation}
\calJ_\eta = \eta\int_{\Omega} H_{2\beta}(\psi) (u-1)^2 \,dx, \quad
\eta>0.
\label{eq:penalty}
\end{equation}
The idea behind is the following. The penalty $\calJ_\eta$ will be
zero where $H_{2\beta}$ is zero or $u=1$ is met.
Near the boundary of $\omega$ the solution $u$ is close to
1, and $\calJ_\eta$ is small. However, in regions like the
one shown in Figure \ref{fig:rsq-1} the value of $\calJ_\eta$ is larger
since $u_h$ is much less than 1.
Notice that instead of the original value used in the state problem, 
 a two times larger value  
$2\beta$ was used in the smoothing of the Heaviside function 
in \eqref{eq:penalty}.

This approach turns out to be effective. In Figures \ref{fig:rsq-2}
and \ref{fig:rsq-2-pot} we show
the results for the same
problem as in Example~\ref{ex:rsq-1} but
with the penalty term included. The value of the penalty parameter was
$\eta=0.1$.  
Now the problematic grey regions near the corners have disappeared, and 
have been replaced by regions where the Heaviside function 
has the value $H_{\beta}=1$. Also the potential $u_h$ is close to 1
near this region (see Figure \ref{fig:rsq-2-pot}).
This time the optimizer needed 1009 iterations and 4098 function evaluations.
The final values are $\calJ = 7.67 \cdot 10^{-5}$ and  $\calJ_\eta =
8.42 \cdot 10^{-5}$.

\begin{figure}[ht!]
\center
\includegraphics[width=0.4\textwidth]{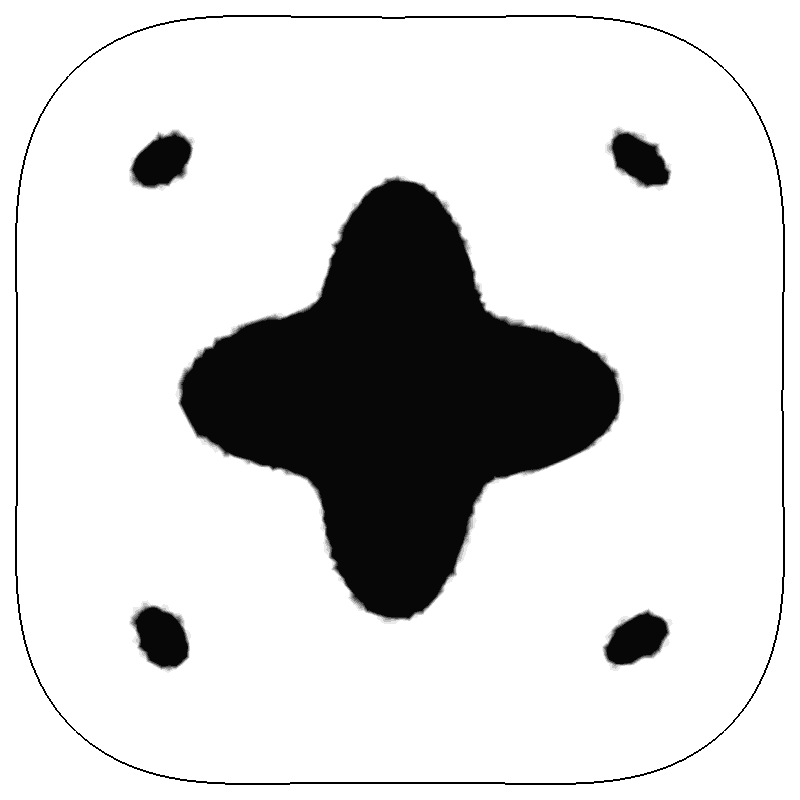}
\caption{The final geometry and $H_\beta(\psi_N^\star)$.}
\label{fig:rsq-2}
\end{figure}

\begin{figure}[ht!]
\center
\includegraphics[width=0.4\textwidth]{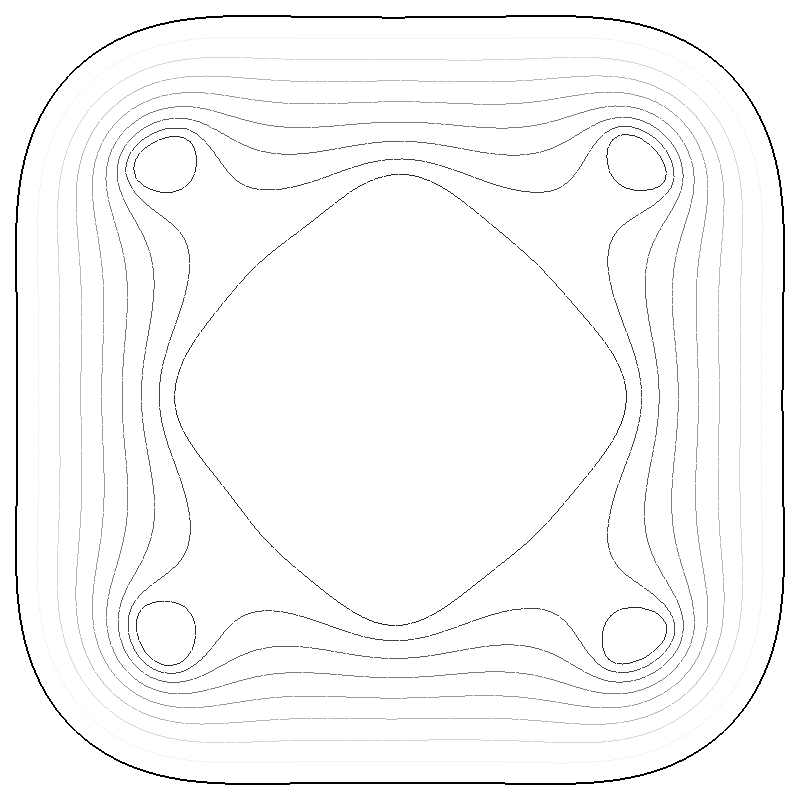}
\caption{The potential $u_h$ related to Example \ref{ex:rsq-2}.}
\label{fig:rsq-2-pot}
\end{figure}



\end{example}

\begin{example} \label{ex:key} \rm
Let $\gamma=-2$ and let the radius function defining the target
boundary $\Gamma_t$  be chosen as 
follows:
$$
   g_t(\theta)=0.5\cos(\theta)+0.8\cos(2\theta)+2, \quad
 \theta\in[0,2\pi[.
$$
The final finite element mesh consists of 15820 linear triangular
elements, and the 
number of RBFs is $15 \times 15$. Value of the penalty parameter was
$\eta = 0.01$. 
The initial value of the objective function (i.e.\ the sum of $\calJ$
and $\calJ_\eta$) was 4.77.
After the total of 356 optimization steps and 1559 function evaluations
this value was reduced to $1.90 \cdot 10^{-5}$.
The results of computations are depicted in Figures~\ref{fig:key} and
\ref{fig:key-pot}. 

This problem was solved in \cite{toi08} using a fixed topology
approach. The results of computations led  to a conclusion that the
inner boundary consists of more 
than one component. From this reason
two holes as an initial approximation of the
inner boundary were introduced ``by hand'', each
parameterized by the radial co-ordinates. The 
obtained result was similar to that in Figure~\ref{fig:key}.

\begin{figure}[ht]
\center
\includegraphics[width=0.7\textwidth]{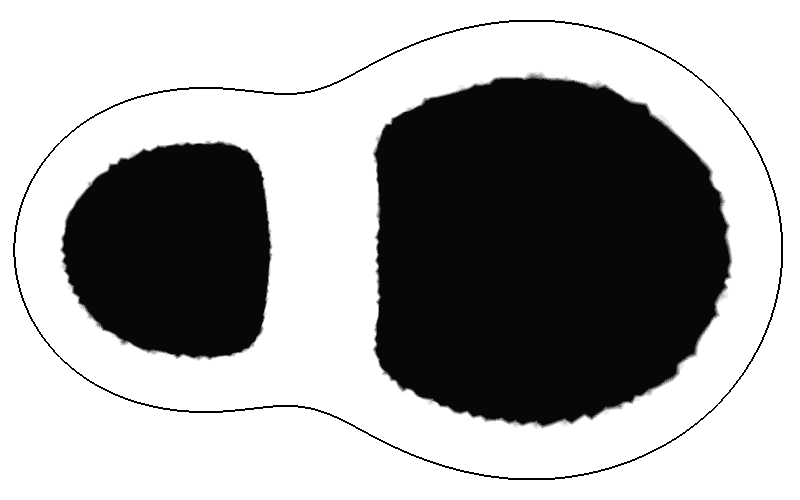} 
\hspace{-0.3truecm}\raisebox{2truecm}{$\Gamma_t$}
\caption{The final geometry and $H_\beta(\psi_N^\star)$}
\label{fig:key}
\end{figure}

\begin{figure}[ht]
\center
\includegraphics[width=0.7\textwidth]{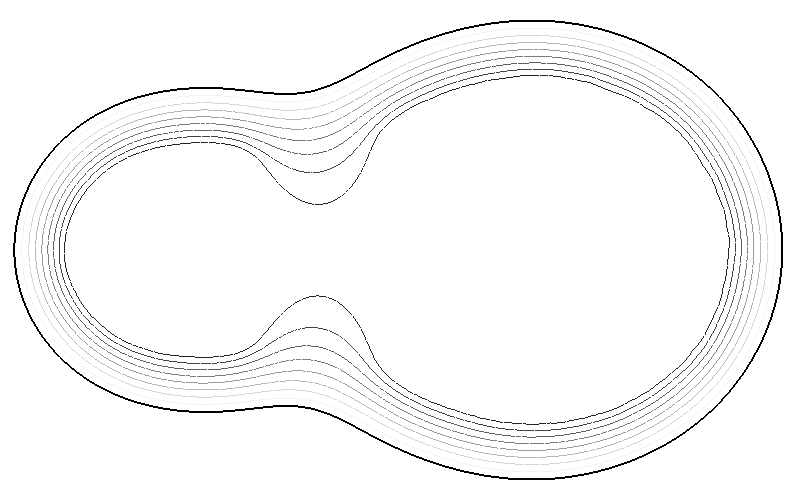}
\caption{The potential $u_h$ related to Example \ref{ex:key}.}
\label{fig:key-pot}
\end{figure}

\end{example}


\section{Conclusions}

In this paper we have considered topology optimization with the state
constraint given by a free boundary problem of Bernoulli type. To
solve efficiently the free boundary problems during the optimization,
the pseudo-solid approach is applied. Its main advantage is that there
is no explicit parametrization of the shape of the free boundary using
e.g.\ splines.
The novelty of the numerical method proposed in this paper is the
combination of the pseudo-solid approach to tackle the free boundary
problem  with a parameterized level set method for shape optimization.
It has been found already in \cite{toi08} that
the problem is very badly conditioned as many different choices of $\omega$ may
lead to nearly identical free boundaries. 
Therefore the progress of the optimization is often
slow. The proposed method can be applied 
in an analogous way to topology optimization problems governed by
other free boundary problems.

\bibliographystyle{ieeetr}
\bibliography{tutkimus}

\section*{Appendix 1}

\subsection*{Smoothing the Heaviside function}

The expression for the $C^2$-smoothed Heaviside function is done by
\begin{equation}
\label{eq:heaviside}
H_\beta(y) =
\begin{cases}
0, &\text{ if } y < -\beta \\
1, &\text{ if } y > \beta \\
\left(\dfrac{y}{\beta} - \dfrac{y^3}{3 \beta^3}\right) \dfrac{3}{4} + \dfrac{1}{2} &\text{ otherwise, }\\
\end{cases}
\end{equation}
where $\beta>0$ is a given
constant.

\begin{figure}[h]
\center
\includegraphics[width=0.8\textwidth]{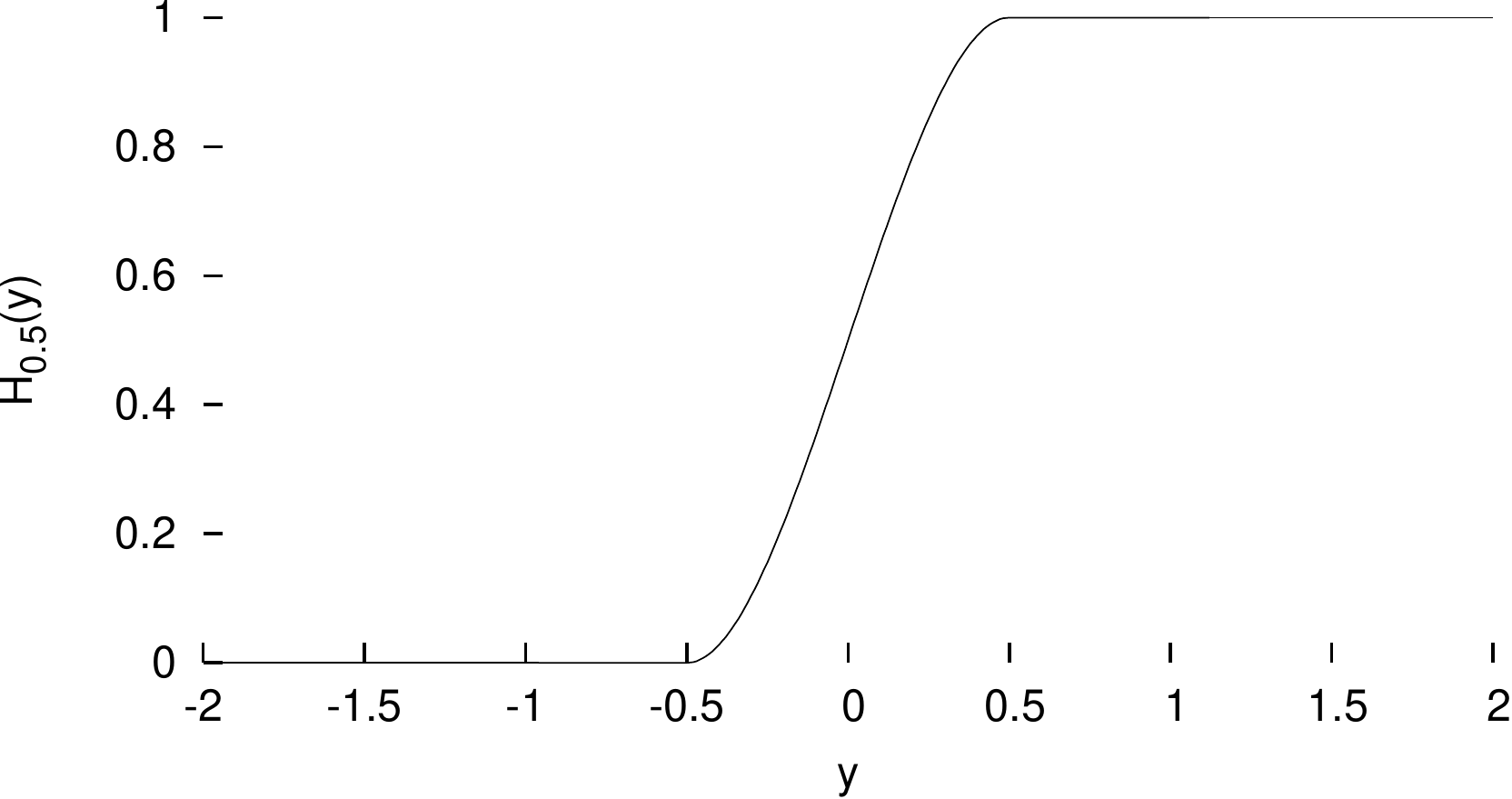}
\caption{The graph of  $H_\beta$, $\beta=\frac12.$} \label{fig:H-beta}
\end{figure}


If $\psi(\xx) \in ]-\beta, \beta[$ then $\xx$ belongs to the so-called
grey region, i.e. $H_\beta(\psi(\xx)) \in ]0,1[$.
Such regions appear near the boundaries of the design domain $\omega$.
To maintain a good quality of the gradients that we obtain, 
we would like to have the geometrical width of the grey region to
be approximately constant everywhere. 
If the scalar function $\psi$ approximated the signed distance function,
this would automatically be the case.
However, we do not make this requirement, but instead we propose the following procedure.

Let $2\delta$ be the desired width of the grey region. Consider
a point $\xx_0$ on the zero level set, and a point $\xx$ located at $\xx = \xx_0 + \delta \nvec$,
where $\nvec$ is the unit normal vector to the zero level set of $\psi$.
The unit normal vector is $\nvec = \nabla \psi / \|\nabla \psi\|$, thus we
can write
\begin{equation}
\xx = \xx_0 + \delta \nvec = \xx_0 + 
  \delta \frac{\nabla \psi (\xx_0)}{\|\nabla \psi(\xx_0)\|}.
\label{eq:h-1}
\end{equation}
Since on the other hand we have the first order approximation 
\begin{equation}
\psi(\xx) \approx \psi(\xx_0) + (\xx-\xx_0) \cdot \nabla \psi(\xx_0),
\label{eq:h-2}
\end{equation}
we get from \eqref{eq:h-1} and \eqref{eq:h-2}  that
\begin{equation}
\psi(\xx)-\psi(\xx_0) \approx  \delta \frac{\nabla \psi(\xx_0)}{\|\nabla
  \psi(\xx_0)\|} \cdot \nabla \psi(\xx_0). 
 = \delta \|\nabla \psi(\xx_0)\|
\end{equation}
Thus if we use
\begin{equation}
\beta(\xx) = \delta \|\nabla \psi(\xx)\|+10^{-6}
\label{eq:beta}
\end{equation}
as the Heaviside parameter,
the width of the grey region will be approximately $2\delta$ everywhere.
The constant $10^{-6}$ is added to prevent division by zero, since
$\beta$ appears in the denominator in \eqref{eq:heaviside}.

\section*{Appendix 2}

\subsection*{Analytical solution in the circular domain}
%

Let $\omega = B(0,R)$, where $B(0, R)$ is the  circle of radius $R$
centered at the origin, 
and $C:=C(R, \gamma)$ be a constant such that
\begin{equation}
\label{eq:CR}
 C \ln(C) - C \ln (R) = -\frac{1}{\gamma}, \quad \gamma<0.
\end{equation}
The function
\begin{equation}
\label{eq:fR}
f(R, \gamma) := C \gamma \ln(\sqrt{x^2+ y^2}) - C \gamma \ln(R) + 1
\end{equation}
satisfies $\Delta f = 0$ in $\Omega \setminus \bomega$,
$f=1$ on $\partial \omega$ and $\nabla f \cdot \nvec = \gamma$
on  $\partial B(0, C(R,\gamma))$.
Thus, given $\gamma<0$ and $\omega = B(0, R)$,
we know that the analytical solution of the free boundary problem
$\eqref{eq:bernoulli}$ 
is $u=f(R, \gamma)$ and $\Omega = B(0, C(R, \gamma))$.
Finally, a level set function representing $B(0,R)$ is given by
\begin{equation}
\label{eq:exscalar}
\hat{\psi}_{R}(x,y) = R - \sqrt{x^2 + y^2}.
\end{equation}
%
%
%
%
%

\end{document}